\documentclass[10pt]{article}

\usepackage[a4paper,textwidth=18cm,textheight=24cm,top=2.85cm, bottom=2.85cm, left=1.5cm, right=1.5cm]{geometry}

\usepackage{template/icdp2009}

\makeatletter
\long\def\@makecaption#1#2{%
\vskip\abovecaptionskip
\sbox\@tempboxa{#1. #2}%
\ifdim \wd\@tempboxa >\hsize
#1. #2\par
\else
\global \@minipagefalse
\hb@xt@\hsize{\box\@tempboxa\hfil}%
\fi
\vskip\belowcaptionskip}
\makeatother

\usepackage{lmodern}

\usepackage{graphicx}
\usepackage{times}
\usepackage{microtype}
\usepackage{graphicx}
\usepackage{subfigure}
\usepackage{booktabs} 
\usepackage{amsmath}
\usepackage{amssymb}
\usepackage[utf8]{inputenc}
\newtheorem{theorem}{Theorem}
\newtheorem{corollary}{Corollary}
\newtheorem{hypothesis}{Hypothesis}

\begin{document}
\noindent

\bibliographystyle{plain}

\title{The Landscape of Multi-Layer Linear Neural Network From the Perspective of Algebraic Geometry}

\authorname{Yang Xiuyi}
\authoraddr{fwinlee@gmail.com}


\maketitle

\begin{abstract}
The clear understanding of the non-convex landscape of neural network is a complex incomplete problem. This paper studies the landscape of linear (residual) network, the simplified version of the nonlinear network. By treating the gradient equations as polynomial equations, we use algebraic geometry tools to solve it over the complex number field, the attained solution can be decomposed into different irreducible complex geometry objects. Then three hypotheses are proposed, involving how to calculate the loss on each irreducible geometry object, the losses of critical points have a certain range and the relationship between the dimension of each irreducible geometry object and strict saddle condition. Finally, numerical algebraic geometry is applied to verify the rationality of these three hypotheses which further clarify the landscape of linear network and the role of residual connection.
\end{abstract}

\section{Introduction}
\label{sec:intro}

The commonly used deep neural network with non-convex loss surface brings significant improvement to many practical applications\cite{krizhevsky2012imagenet}. 
The difficulty of non-convex optimization was manifest in the practical development of early neural networks \cite{blum1992training}. In the past few years, with the introduction of some new structures, such as residual connection \cite{he2016deep}, normalization techniques \cite{ioffe2015batch,ba2016layer} and so on, non-convex neural networks optimized by stochastic gradient descent(SGD) and its variants can often get very low loss value in practice. The different non-convex landscapes caused by different components are the key to understand why SGD and its variants works.

It has been demonstrated that the linear network is similar to the nonlinear network in many aspects. The work of \cite{saxe2013exact} analyzed the learning dynamics of linear network and non-linear network have similar patterns. Components of nonlinear networks such as residual connection and normalization method can also be used to solve the gradient vanishing or explosion problem of linear networks. These similarities are the first reason why we regard linear network as an ideal model for the research of non-linear neural network. Secondly, general results of nonlinear neural network are often difficult to derive analytically, whereas the linear network can often be analyzed in mathematical detail \cite{baldi1995learning}. So the main topic in present paper is exploring the loss surface of linear (residual) network.

These results that every local minimum is a global minimum and every critical point that is not a global minimum is a saddle point are attained under different assumptions \cite{kawaguchi2016deep,lu2017depth,laurent2018deep,hardt2016identity}. The works of \cite{yun2017global,zhou2017critical} present conditions for a critical point of the risk function to be a global minimum. These results demonstrated that saddle points of linear network are the biggest barrier for SGD. We will use algebraic geometry to zoom in the critical points of linear network to explain the following experimental phenomena. 

Suppose input and output data are fitted perfectly by a linear network, the parameters are initialized with random Gaussian with zero means, SGD is used to minimize squared error risk. When the depth of neural network is not so large, the loss of finial solution of optimization is close to zero. When the gradient disappears due to the increase of depth, the residual connection can make the loss of the final solution of linear network close to zero. Three questions arise here. One is what structure near point $\boldsymbol{0}$ causes the gradient to disappear, another is how the residual connection resolves the dying gradient issue, and the final one is what structure the saddle points have to explain the saddle points is not an obstacle to SGD?

\section{Linear (residue) network}
\label{sec:linear network}

In this section, we describe notations for multi-layer linear neural networks and multi-layer linear residue networks and the problem formulation.

Let's say that there are $m$ training data, $x_i, y_i $ is the ith input data and the corresponding output. $X \in \mathbb{R}^{d_x \times m}$  be the data matrix and $Y \in \mathbb{R}^{d_y \times m}$  be the target matrix, $d_x,d_y$ are the number of input and output units respectively. 
Since the multi-layer linear network (number of layers greater than 1) is a forward neural network with identity mapping activation function, if do not consider the residue connection, normalization methods, etc., and use square error loss, the objective function to be optimized is

\begin{equation}
\label{eq1}
L(W) = \mathrm{\frac{1}{2} \sum_{i=1}^{m} \Vert \prod_{k=1}^{H+1} W_k x_i - y_i \Vert_F^2}
\end{equation}

Where $\prod_{k=1}^{H+1} W_k = W_{H+1} W_{H} \cdots W_1, n>1$, for $k=1, ..., H+1$, $W_k \in \mathbb{R}^{d_k \times d_{k-1}}$ is the weights between adjacent layers. For notiaonal simplicity we
let $d_0 = d_x$ and $d_{H+1} = d_y$ and the width of net is defined as $k= \text{min}(d_0, d_1, \cdots, d_{H+1})$. The number of weights, or variables, $n=d_0 \times d_1 + d_1 \times d_2 + \cdots + d_H \times d_{H+1} $.

For linear residual networks, we only consider the case that the units of each layer are equal and residue connection only skip one layer as in\cite{hardt2016identity}. The empirical risk in this case is as follows,

\begin{equation}
\label{eq2}
L_r(W^{\prime}) = \mathrm{\frac{1}{2} \sum_{i=1}^{m} \Vert \prod_{k=1}^{H+1} (I+W_k^{\prime}) x_i - y_i \Vert_F^2}
\end{equation}

It is easy to conclude that the loss surface of linear residual network is the translation of the corresponding linear network because that the objective function equation (\ref{eq2}) is obtained by reparameterization $W_k = I+W_k^{\prime}$ from equation (\ref{eq1}). Therefore, this conclusion that there are only global minima and saddle points holds for linear residual networks.

Let the partial derivative of $L(W)$ with respect to $W_k$ equals to $0$,

\begin{equation}
\label{eq3}
\begin{aligned}
 \frac{\partial L(W)}{\partial W_i} &= ( \prod_{k=i+1}^{H+1} W_k ) ^T (\prod_{k=1}^{H+1} W_k X - Y) X^T ( \prod_{k=1}^{i-1} W_k ) ^T \\ &= 0    
\end{aligned}
\end{equation}

$ \text{ for \ } i = 1, \ldots , H+1.$

Thus, this is a set of polynomial equations in the entries
of $W_1,...,W_{H+1}$, so is linear residual networks. \cite{mehta2018loss} remove all the flat stationary points by adding an extension of $L_2$-regularization to the loss function (\ref{eq1}), then pick out real solutions from all isolated complex points and find that there are indeed local minima which are not global minima contrary to the available conclusions in the unregularized case \cite{kawaguchi2016deep}. As proved in \cite{taghvaei2017regularization}, regularization alters the loss surface of linear network.

Not only the above gradient equations (\ref{eq3}) can be transformed into polynomial equations, but also some other results about linear networks can be transformed into polynomial equations. A recent major work on linear networks\cite{yun2017global}, gives the condition that the critical points are globally optimal when the width of linear network equals to the number of neuron of input or output layer. This condition can be considered from the perspective of polynomial equations. 

\begin{theorem}
\label{thm:global condition}
\cite{yun2017global} If $k = \text{min} \{ d_x, d_y \} $, define the following set
$V_1 := \{ (W_1,..., W_{H+1}): \text{rank}(W_{H+1} \cdots W_1) = k \}$.
Then, every critical point of $L(W)$ in $V_1$ is a global minimum. Moreover, every critical point of
$L(W)$ in $V^c_1$
is a saddle point.
\end{theorem}


Here the solution set of $\text{rank}(W_{H+1} \cdots W_1) = k$ is  complementary set of solution set of $\text{rank}(W_{H+1} \cdots W_1) < k$. The latter is equivalent to the determinants of all the $k \times k$ submatrix of $W_{H+1} \cdots W_1$ equal to zero. So saddle points are solutions of polynomial equations. When solutions of gradient equations and these polynomial equations lives in $\mathbb{C}^n$, they are called \textbf{complex saddle points}. The corresponding global minima are called \textbf{complex global minima}. These saddle points are further divided into different subsets each of which satisfies $\text{rank}(W_n \cdots W_1) = i, i=0,...,k-1$. Likewise, these different subsets can be computed through $\text{rank}(W_n \cdots W_1) < i, i=1,...,k$. There must have a trivial solution that is saddle point 0 when $\text{rank}(W_n \cdots W_1) = 0$ and the losses at these saddle points of rank zero equal to the loss at saddle point 0.

Since global minima are critical points which meet the condition $\text{rank}(W_{H+1} \cdots W_1) = k$, we have following hypothesis,

\begin{hypothesis}
\label{conj:1}

When the width of nets equals to the number of input or output units, the larger the rank of the product of the chained product of weight matrices, the smaller the loss of the corresponding critical points.

\end{hypothesis}

The critical points under condition $\text{rank}(W_n \cdots W_1) = i$ may include many stationary points or none. So the comparison of the loss of critical points which satisfy $\text{rank}(W_n \cdots W_1) = i$ and $\text{rank}(W_n \cdots W_1) = i-1$ is based on the rule: the loss of any point in one set is greater than or equal to or less than or equal to that of any point in the other set. In the case of one of them is empty set, the comparison is trivially true.

\section{Methodology}
\label{set:method}
A brief introduction of algebraic geometry is provided in subsection \ref{sec:AG}. Then in next subsection \ref{sec:relation} the reason why complex field is used to solve equations of critical points is explained. In subsection \ref{sec:1t1u}, a very simple network that is known all the stationary points is used as an example to explain the algebraic geometry terms in present paper. This example leads to two other hypotheses related to the properties of critical points. Finally in subsection \ref{sec:NAG}, three methods in numerical algebraic geometry are introduced, which will be used to compute complex critical points.

\subsection{A brief introduction to algebraic geometry}
\label{sec:AG}
The fundamental goal of algebraic geometry is to study solution sets of systems of polynomial equations in several variables. 
For a subset $S \subset \mathbb{C}[x_1,\cdots, x_n]$ of polynomials we define

$V(S) := \{ x \in \mathbb{C}^n: f(x) = 0 \textit{ for all } f \in S \}$

\textbf{closed affine varieties} of $\mathbb{C}^n$.

In this paper, the closed affine varieties to be studied are, the complex solutions of gradient equations of linear (residue) network(\textbf{complex critical points}) and its subsets. In view of the importance of the critical point in the study of the loss surface of linear network, we name complex critical points as \textbf {linear neural variety}. The word closed appears in the definition because it can be shown that the closed affine varieties of $\mathbb{C}^n$ satisfy the axioms to be the closed subsets for a topology on $\mathbb{C}^n$. This topology is called the Zariski topology on $\mathbb{C}^n$. The induced topology on a closed affine variety $V$ of $\mathbb{C}^n$ is called the \textbf{Zariski topology} on $V$. The Zariski topology is coarser than the usual complex topology. Irregular (non-smooth or local dimension greater than the dimension of variety) points are measure 0 in variety $V$ is deduced from that smooth points on this variety $V$ are zariski dense in it. A \textbf{generic property} is a property which is true for almost every point of a variety and a \textbf{generic point} of a variety is a point at which all generic properties are true. For example, a generic point of a variety is a smooth point.

One of main goals in present paper is decomposing a variety into finite union of irreducible varieties. Let $V$ be a topological space. We say that $V$ is \textbf{reducible} if it can be written as $V = V_1 \cup V_2$ for closed subsets $V_1, V_2 \subsetneq V$. Otherwise $V$ is called \textbf{irreducible}. The finite irreducible decomposition of an affine variety is $V = V_1 \cup V_2 \cdots \cup V_r$ of irreducible closed subsets, up to permutation. $V_i,i=1,...,r$ are called the irreducible components of $V$. These irreducible closed subsets can be written as the zero locus of finitely many polynomials by Hilbert’s Basis Theorem. If the components of linear neural variety belong to complex global minima they are called \textbf{global minimum components}, if belongs to complex saddle points, then are called \textbf{saddle components}. Each algebraic component has a well-defined dimension and degree. Every irreducible dimension 0 algebraic component is a set of single point. A irreducible algebraic curve has dimension 1 and so on. The \textbf{degree} of an affine variety of dimension $n$ is the number of intersection points of the variety with $n$ generic hyperplanes. The number of intersection is counted with intersection multiplicity. Complex irreducible component $V_i$ has a important property that if $V_i$ is an irreducible algebraic variety over $\mathbb{C}$, then $V_i(\mathbb{C})$ is connected.

\subsection{Why linear neural variety lives in $\mathbb{C}^{n}$}
\label{sec:relation}
From the definition of linear neural variety, it can be seen that real stationary points are the intersection of linear neural variety and $\mathbb{R}^n$, that is real solutions of gradient system. Symbolic computational methods such as the Gröbner basis \cite{cox2013ideals,cox2006using} and cylindrical algebraic decomposition \cite{jirstrand1995cylindrical} in semi-algebraic geometry can be used to solve this system by identifying it as polynomial equations. But Gröbner basis and cylindrical algebraic decomposition may have a worst case complexity doubly exponential in the number of solutions of the polynomial system and the size of the input,respectively. Due to algorithmic complexity issues of these two methods, \cite{mehta2018loss} employ the numerical homotopy continuation method \cite{sommese2005numerical,bates2013numerically} to sort these purely real solutions out from the complex solutions through removing flat stationary points by a generalized $L_2$ penalty term. However, the loss surface is changed after the introducing of this penalty term. Therefore, in present paper, we consider the linear neural variety over the complex field and studying the properties of critical points by studying the corresponding linear neural variety.

In particular, there is a theorem in algebraic geometry relates complex solutions of polynomial equations with real solutions. That is,

\begin{theorem}
\label{thm:sottile}
\cite{sottile2016real} Let $V_i \subset \mathbb{C}^n$ be an irreducible variety defined by real polynomials. If $V$ has a smooth real point, then $V_i(\mathbb{R}) = V_i \cap \mathbb{R}^n$ is Zariski dense in $V$. 
\end{theorem}


Training data in present paper are over the real number, so if linear neural variety can be decomposed into finite irreducible components, each of which are defined by finite real polynomials and if each irreducible component has a smooth real point, then all algebraic and geometric information about this component is already captured in real critical points, and vice-versa. So, if the complex dimension of each irreducible component are known, then corresponding real counterpart is with the same real dimension. At the same time, each irreducible component is regarded as an individual just like the isolated critical point, its pseudo loss and pseudo eigenvalues of Hessian matrix should be studied. The word pseudo appears here due to that each irreducible closed variety is solved over complex number, yet loss and eigenvalues of Hessian are only meaningful for real point of a component at first glance. It will be seen that loss and eigenvalues of real critical points $V_i(\mathbb{R})$ can be reflected by corresponding irreducible component $V_i$. After all, both complex critical points and real critical points correspond to the same polynomial equations, only one are solved over complex field, the other are solved over real number.

\subsection{Example: a very simple linear neural network}
\label{sec:1t1u}

This section will first examine the linear variety of simple network with width 1, then the loss and eigenvalues of Hessian matrix of each irreducible component are derived analytically, this leads to two other main hypotheses in present paper.

the gradient polynomial system with respect to the loss function is:
\begin{equation}\label{eq4}
\frac{\partial l(w)}{\partial w_{i}} =  -(y - (\prod_{k=1}^{H+1} w_k) x)((\prod_{k=1}^{i-1} w_k) x \times \prod_{k=i+1}^{H+1} w_k) = 0
\end{equation}

The quadratic loss is globally minimized by solutions of the equation $(\prod_{k=1}^{H+1} w_k) x - y = 0$. On the contrary, if $(\prod_{k=1}^{H+1} w_k x - y) \neq 0$, then the locus of 

\begin{equation}\label{eq5}
\prod_{k=1}^{i-1} w_k x \times \prod_{k=i+1}^{H+1} w_k = 0, i =1,...,n;
\end{equation}

be composed of 

\begin{equation}\label{eq6}
w_i = 0,
w_j = 0, i,j \in \{ 1,...,H+1 \}, i \neq j;
\end{equation}

With terms in algebraic geometry, above is equivalent to the linear variety defined by equation (\ref{eq4}) is decomposed as the union of global minimum variety and $\binom{H+1}{2}$ saddle components defined by equations (\ref{eq5}) and (\ref{eq6}) respectively. The global minimum component is a complex surface with complex dimension $H=n-1$, each irreducible component is a complex plane with complex dimension $H-1=n-2$ and these saddle component intersect at saddle point $\boldsymbol{0}$. Note that the real dimension is consistent with the complex dimension, these components can be treated as a real surface and a real plane with real dimension $H=n-1$ and $H-1=n-2$ respectively.

\begin{figure}[ht]
\vskip 0.2in
\begin{center}
\centerline{\includegraphics[width=\columnwidth]{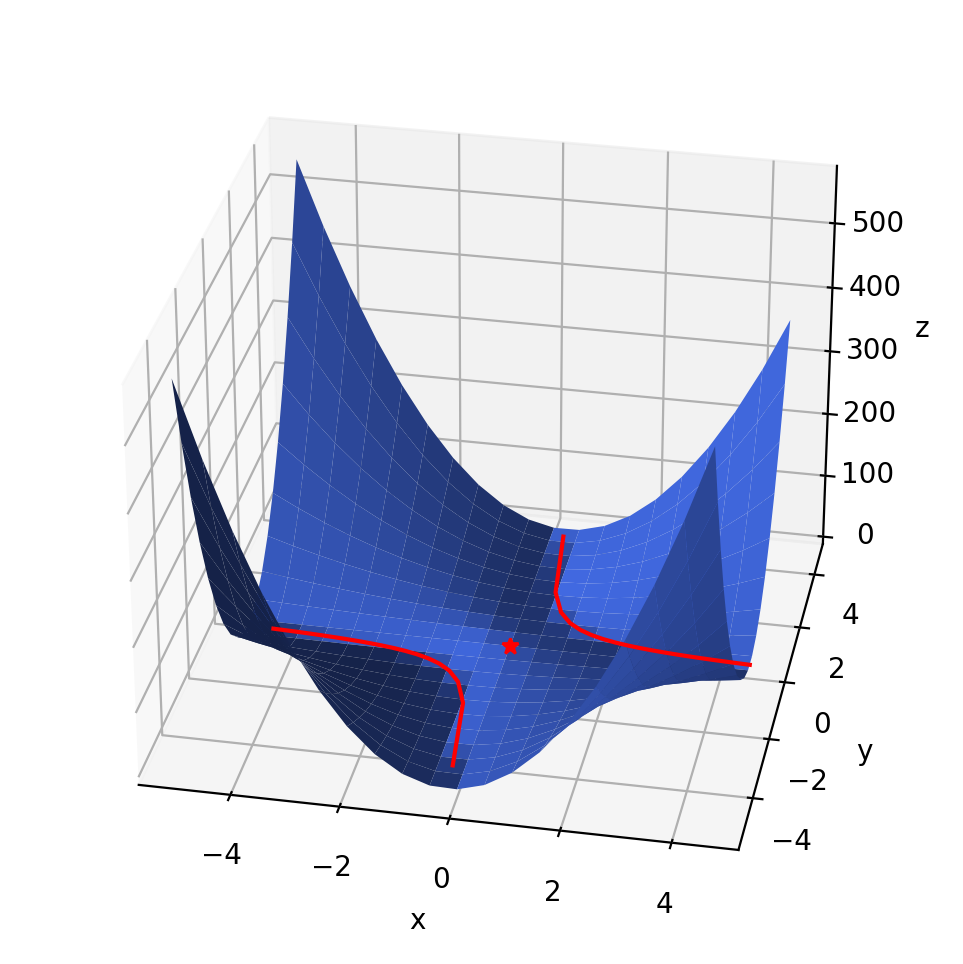}}
\caption{The loss surface of $z =(1-x y)^2$, the set of global minima is curve $x y=1$, the only saddle point is 0.}
\label{fig:1}
\end{center}
\vskip -0.2in
\end{figure}

Observe that the pseudo loss of each irreducible component is constant at every complex point, this pseudo loss is the same as the actual loss at real critical point. So, there are only two loss values, one is global minimum 0 at global minimum component, the other is saddle loss $\mathrm{\frac{1}{2} \sum_{i=1}^{m} (y_i)^2}$ at each saddle component. The loss surface of 2 units neural network is seen as in figure \ref{fig:1}

Now, it is time to compute the eigenvalues of each irreducible components.
First, the elements of Hessian matrix is 
\begin{equation}\label{eq7}
\left
\{
\begin{array}{lr} 
\begin{aligned}
\frac{\partial l^2(w)}{\partial w_{i}^2} &= x^2 (\prod_{k\neq i}  w_k)^2 \\
\frac{\partial l^2(w)}{\partial w_i w_j} &=  -x y\prod_{k \neq i,j} w_k+2 x^2 \prod w_k \prod_{k\neq i,j}  w_k

\end{aligned}
\end{array} 
\right.
\end{equation}

The only non constant zero eigenvalues of Hessian at a point on global minimum component is $\lambda =\sum_k  (y^2/w_k^2)$ and the two non constant zero eigenvalues of Hessian at a point on saddle components are $\lambda = \pm (-x y\prod_{k \neq i,j} w_k+2 x^2 \prod w_k \prod_{k\neq i,j}  w_k)$. Around saddle point $\boldsymbol{0}$, there is a geometrically flat zone caused by small eigenvalues due to the eigenvalues of point $\boldsymbol{0}$ is all zero and two non constant zero eigenvalues are continuous function of points on saddle components. 
Whether real critical point or not on a irreducible component satisfy the same eigenvalue function. Another point to note is that the sum of the dimension of a irreducible component and the number of non constant zero eigenvalues equals to $n=H+1$.

Two hypotheses are proposed according two properties of this simple network that are constant loss on a irreducible component and the relation between dimension and the number of constant zero eigenvalues of a irreducible component.

\begin{hypothesis}
\label{conj:constant loss}
The loss function of multi-layer linear network is constant on each irreducible algebraic component of its linear neural variety.
\end{hypothesis}

This hypothesis has a very obvious corollary,
\begin{corollary}
\label{cor:1}
The loss function of multi-layer linear network is constant on finite many connected algebraic components of its linear neural variety.
\end{corollary}

The hypothesis (\ref{conj:1}) make it possible to compute the loss value of every irreducible component, that is the loss of every critical points as long as we sample a complex or real point of this component. The meaning of the corollary (\ref{cor:1}) is that saddle point $\boldsymbol{0}$ in many nets is the intersection point of several irreducible components, this will be seen in section \ref{sec:computation}, so the losses of these several components will be equal to the loss of point $\boldsymbol{0}$.

\begin{hypothesis}
\label{conj:eigenvalue}
For a generic point of each irreducible component, the number of constant zero eigenvalues of Hessian at this point equals to the dimension of this component.
\end{hypothesis}

This hypothesis is equivalent to say that the characteristic polynomial $f(W, \lambda)$ of Hessian at this component has a factor $\lambda^r$, $f(W, \lambda)= \lambda^r g(W,\lambda)$, $r$ is the dimension of this component. The coefficients of $f$ and $g$ are (rational) polynomials of parameters of network. If some point of this component makes the constant term of $g$ vanish, then the number of constant zero eigenvalues is $r+1$, etc. These points is rare due to they are the closed subvariety of this component. When several components meet at a point, the number of constant zero eigenvalues may increase due to different component correspond to different eigenvalues of Hessian matrix as above network example. So, we have the following corollary,

\begin{corollary}
\label{cor:2}
For every point of each irreducible component, the number of constant zero eigenvalues of Hessian is greater than or equals to the dimension of this component. For the intersection point of many irreducible components, the number of constant zero eigenvalues at intersection point is greater than or equals to the dimension of any of these components.
\end{corollary}

Hypothesis (\ref{conj:eigenvalue}) and its corollary (\ref{cor:2}) means that we can sample a random point of a component to compute its eigenvalues so that we can get the number of constant zero eigenvalues, this number must be consistent with the dimension of the component. Due to $\boldsymbol{0}$ in many cases is the intersection point, there may have more constant zero eigenvalues at point $\boldsymbol{0}$, so the area around 0 of loss surface may be flatter than other areas.

\subsection{A brief introduction to numerical algebraic geometry}
\label{sec:NAG}
Numerical algebraic geometry\cite{sommese2005numerical,bates2013numerically} is a subject which uses methods from numerical analysis to manipulate varieties. The primary computational method used in numerical algebraic geometry is the numerical homotopy continuation, which solve a system of polynomial equations from the known solution of another system of polynomial equations. In order to manipulate high-dimensional solution of polynomial equations, the data structure "witness set" is been introduced to encode information for algebraic varieties. The witness set for a pure dimensional variety(the dimension of all irreducible component is equal) $V$ contains three part$(S, L, S \cap L)$. The first part $S$ is a system of polynomial equations which define the studied variety, the second part $L$ is a generic slicing plane the dimension of which is codimension of $V$, the last part $S \cap L$ is witness points(the intersection points of slicing plane $L$ and variety $V$) the number of which equals to the degree of variety $V$. 

Witness sets are been used to compute anything of interest in numerical algebraic geometry such as the dimension and degree of a variety, the intersection of two varieties, etc. In present paper, we are interested in three methods using witness sets: numerical irreducible decomposition, component membership testing, and component sampling. \textbf{numerical irreducible decomposition} is computing a witness set for each of the irreducible components, \textbf{component membership testing} method tests whether a given point is on a variety and \textbf{component sampling} method is used to sample a generic point on a variety. These three methods are been implemented in $Bertini$ \cite{bates2013numerically}.

With these three methods of numerical algebraic geometry, we can compute the loss and the number of constant zero eigenvalues of a irreducible component,

\begin{enumerate}
\item First represent linear neural variety as the union of irreducible components by numerical irreducible decomposition method.

\item Then sample a generic point on a component to compute the loss and eigenvalues at this point, this will lead to results we want.

\end{enumerate}

\section{Experiments}
\label{sec:computation}
In this section, we will compute linear neural varieties of some networks, then compute loss on each irreducible component. In order to verify the rational of the hypothesis \ref{conj:constant loss} and its corollary \ref{cor:1}, we sample at least two points to compute their losses to get the final loss of a component.  We use component membership testing to test point $\boldsymbol{0}$ is actually the intersection point of several components in some cases of nets and the losses of these intersecting components indeed have the same loss. Then the rational of the hypothesis \ref{conj:1} is obvious by solving the saddle points under different ranks. The checking of reasonableness of hypothesis \ref{conj:eigenvalue} and its corollary \ref{cor:2} is also obvious after we computed each irreducible component. Finally, the validity of the loss surface of linear residue network is the translation of corresponding linear network by checking whether two linear neural varieties are equal after translate one linear neural variety, component membership testing method will be used during checking the relationship between two sets.

Experimental settings are as follows. Due to we have limited computational resource, we only compute width and depth of nets less than Four. Algebraic geometry method is suitable for linear network trained with any data, so we don't make any assumptions about the amount of data or the distribution of data. But in order to check whether the loss value of global minimum component is zero, the input and output data can be perfectly fitted by linear network. For convenience and the loss surface is non-convex regardless of the number of training data, the training data is up to three. 

\subsection*{Linear neural varieties and losses}
The results of numerical irreducible decomposition of linear neural varieties and losses are summarized in Table \ref{table:1-1-1-1-1-critical}, \ref{table:2-2-2-1-critical},  \ref{table:2-2-2-1-critical-2data}, \ref{table:2-2-2-1-critical-3data},  \ref{table:2-2-2-critical-2data}, \ref{table:2-2-3-critical-2data}, \ref{table:2-1-2-2-critical-2data}, \ref{table:3-2-3-critical-2data}, \ref{table:2-1-2-critical-2data}. We observe that in all of these results, there is only one global minimum component and the losses of components which saddle point $\boldsymbol{0}$ is on are larger than that of any other component.

\begin{table*}
\centering
\begin{minipage}[t]{\columnwidth}
\caption{The linear neural variety of 1-1-1-1-1 linear net, with 2 training data}
\label{table:1-1-1-1-1-critical}
\vskip 0.15in
\begin{center}
\begin{tabular}{|c|c|c|ll}
\cline{1-3}
\multicolumn{3}{|c|}{1-1-1-1-1, 2 data}                                                                      &  &  \\ \cline{1-3}
Irreducible components  & Loss & \begin{tabular}[c]{@{}c@{}}Is saddle $\boldsymbol{0}$ on this \\ component?\end{tabular} &  &  \\ \cline{1-3}
Dim 3, deg 4, 1 component & 0    & No                                                                        &  &  \\ \cline{1-3}
Dim 2, deg 1, 6 components & 10   & Yes                                                                       &  &  \\ \cline{1-3}
\end{tabular}
\end{center}
\vskip -0.1in
\end{minipage}\hfill
\begin{minipage}[t]{\columnwidth}
\caption{The linear neural variety of 2-2-2-1 linear net, with 1 training data}
\label{table:2-2-2-1-critical}
\vskip 0.15in
\begin{center}
\begin{tabular}{|c|c|c|ll}
\cline{1-3}
\multicolumn{3}{|c|}{2-2-2-1, 1 data}                                                                       &  &  \\ \cline{1-3}
Irreducible components   & Loss & \begin{tabular}[c]{@{}c@{}}Is saddle $\boldsymbol{0}$ on this \\ component?\end{tabular} &  &  \\ \cline{1-3}
Dim 9, deg 3, 1 component & 0    & No                                                                        &  &  \\ \cline{1-3}
Dim 6, deg 1, 1 component  & 12.5 & Yes                                                                       &  &  \\ \cline{1-3}
Dim 6, deg 3, 2 components  & 12.5 & Yes                                                                       &  &  \\ \cline{1-3}
\end{tabular}

\end{center}
\vskip -0.1in
\end{minipage}

\bigskip

\begin{minipage}[t]{\columnwidth}
\caption{The linear neural variety of 2-2-2-1 linear net, with 2 training data}
\label{table:2-2-2-1-critical-2data}
\vskip 0.15in
\begin{center}
\begin{tabular}{|c|c|c|ll}
\cline{1-3}
\multicolumn{3}{|c|}{2-2-2-1, 2 data}                                                                     &  &  \\ \cline{1-3}
Irreducible components & Loss & \begin{tabular}[c]{@{}c@{}}Is saddle 0 on this \\ component?\end{tabular} &  &  \\ \cline{1-3}
Dim 8, deg 9, 1 component      & 0    & No                                                                        &  &  \\ \cline{1-3}
Dim 6, deg 1, 1 component      & 13   & Yes                                                                       &  &  \\ \cline{1-3}
Dim 6, deg 3, 1 component      & 13   & Yes                                                                       &  &  \\ \cline{1-3}
Dim 6, deg 3, 1 component      & 13   & Yes                                                                       &  &  \\ \cline{1-3}
\end{tabular}
\end{center}
\vskip -0.1in
\end{minipage}\hfill
\begin{minipage}[t]{\columnwidth}
\caption{The linear neural variety of 2-2-2-1 linear net, with 3 training data}
\label{table:2-2-2-1-critical-3data}
\vskip 0.15in
\begin{center}
\begin{tabular}{|c|c|c|ll}
\cline{1-3}
\multicolumn{3}{|c|}{2-2-2-1, 3 data}                                                                     &  &  \\ \cline{1-3}
Irreducible components & Loss & \begin{tabular}[c]{@{}c@{}}Is saddle $\boldsymbol{0}$ on this \\ component?\end{tabular} &  &  \\ \cline{1-3}
Dim 8, deg 9, 1 component      & 0    & No                                                                        &  &  \\ \cline{1-3}
Dim 6, deg 1, 1 component     & 13.5 & Yes                                                                       &  &  \\ \cline{1-3}
Dim 6, deg 3, 1 component      & 13.5 & Yes                                                                       &  &  \\ \cline{1-3}
Dim 6, deg 3, 1 component      & 13.5 & Yes                                                                       &  &  \\ \cline{1-3}
\end{tabular}
\end{center}
\vskip -0.1in
\end{minipage}

\bigskip

\begin{minipage}[t]{\columnwidth}
\caption{The linear neural variety of 2-2-2 linear net, with 2 training data}
\label{table:2-2-2-critical-2data}
\vskip 0.15in
\begin{center}
\begin{tabular}{|c|c|c|ll}
\cline{1-3}
\multicolumn{3}{|c|}{2-2-2, 2 data}                                                                       &  &  \\ \cline{1-3}
Irreducible components & Loss & \begin{tabular}[c]{@{}c@{}}Is saddle $\boldsymbol{0}$ on this \\ component?\end{tabular} &  &  \\ \cline{1-3}
Dim 4, deg 8, 1 component      & 0    & No                                                                        &  &  \\ \cline{1-3}
Dim 3, deg 2, 1 component      & 5.3  & No                                                                        &  &  \\ \cline{1-3}
Dim 3, deg 2, 1 component      & 1.69 & No                                                                        &  &  \\ \cline{1-3}
Dim 0, deg 1, 1 component      & 7    & Yes                                                                       &  &  \\ \cline{1-3}
\end{tabular}
\end{center}
\vskip -0.1in
\end{minipage}\hfill
\begin{minipage}[t]{\columnwidth}
\caption{The linear neural variety of 2-2-3 linear net, with 2 training data}
\label{table:2-2-3-critical-2data}
\vskip 0.15in
\begin{center}
\begin{tabular}{|c|c|c|ll}
\cline{1-3}
\multicolumn{3}{|c|}{2-2-3, 2 data}                                                                       &  &  \\ \cline{1-3}
Irreducible components & Loss & \begin{tabular}[c]{@{}c@{}}Is saddle $\boldsymbol{0}$ on this \\ component?\end{tabular} &  &  \\ \cline{1-3}
Dim 4, deg 8, 1 component        & 0    & No                                                                        &  &  \\ \cline{1-3}
Dim 4, deg 4, 1 component        & 2    & No                                                                        &  &  \\ \cline{1-3}
Dim 4, deg 4, 1 component         & 7    & No                                                                        &  &  \\ \cline{1-3}
Dim 2, deg 1, 1 component      & 9    & Yes                                                                       &  &  \\ \cline{1-3}
\end{tabular}
\end{center}
\vskip -0.1in
\end{minipage}

\bigskip

\begin{minipage}[t]{\columnwidth}
\caption{The linear neural variety of 2-1-2-2 linear net, with 2 training data}
\label{table:2-1-2-2-critical-2data}
\vskip 0.15in
\begin{center}
\begin{tabular}{|c|c|c|ll}
\cline{1-3}
\multicolumn{3}{|c|}{2-1-2-2, 2 data}                                                                     &  &  \\ \cline{1-3}
Irreducible components & Loss & \begin{tabular}[c]{@{}c@{}}Is saddle $\boldsymbol{0}$ on this \\ component?\end{tabular} &  &  \\ \cline{1-3}
Dim 5, deg 1, 1 component      & 13   & Yes                                                                       &  &  \\ \cline{1-3}
Dim 5, deg 2, 1 component      & 13   & Yes                                                                       &  &  \\ \cline{1-3}
Dim 5, deg 3, 1 component      & 13   & Yes                                                                       &  &  \\ \cline{1-3}
Dim 5, deg 6, 1 component      & 0    & No                                                                        &  &  \\ \cline{1-3}
Dim 4, deg 1, 1 component     & 13   & Yes                                                                       &  &  \\ \cline{1-3}
\end{tabular}
\end{center}
\vskip -0.1in
\end{minipage}\hfill
\begin{minipage}[t]{\columnwidth}
\caption{The linear neural variety of 3-2-3 linear net, with 2 training data}
\label{table:3-2-3-critical-2data}
\vskip 0.15in
\begin{center}
\begin{tabular}{|c|c|c|ll}
\cline{1-3}
\multicolumn{3}{|c|}{3-2-3, 2 data}                                                                        &  &  \\ \cline{1-3}
Irreducible components & Loss  & \begin{tabular}[c]{@{}c@{}}Is saddle $\boldsymbol{0}$ on this \\ component?\end{tabular} &  &  \\ \cline{1-3}
Dim 6, deg 4, 1 component      & 8.17  & No                                                                        &  &  \\ \cline{1-3}
Dim 6, deg 4, 1 component      & 97.82 & No                                                                        &  &  \\ \cline{1-3}
Dim 6, deg 8, 1 component      & 0     & No                                                                        &  &  \\ \cline{1-3}
Dim 4, deg 1, 1 component      & 106   & Yes                                                                       &  &  \\ \cline{1-3}
\end{tabular}
\end{center}
\vskip -0.1in
\end{minipage}

\bigskip

\begin{minipage}[t]{\columnwidth}
\caption{The linear neural variety of 2-1-2 linear net, with 2 training data}
\label{table:2-1-2-critical-2data}
\vskip 0.15in
\begin{center}
\begin{tabular}{|c|c|c|ll}
\cline{1-3}
\multicolumn{3}{|c|}{2-1-2, 2 data}                                                                       &  &  \\ \cline{1-3}
Irreducible components & Loss & \begin{tabular}[c]{@{}c@{}}Is saddle $\boldsymbol{0}$ on this \\ component?\end{tabular} &  &  \\ \cline{1-3}
Dim 1, deg 1, 2 components & 13   & Yes                                                                       &  &  \\ \cline{1-3}
Dim 1, deg 2, 1 component      & 0    & No                                                                        &  &  \\ \cline{1-3}
\end{tabular}
\end{center}
\vskip -0.1in
\end{minipage}

\end{table*}

\begin{table*}
\centering

\begin{minipage}[t]{\columnwidth}
\caption{The saddle components of 2-2-2-1 linear net, with 2 training data}
\label{table:saddle-2-2-2-1-critical-2data}
\vskip 0.15in
\begin{center}
\begin{tabular}{|c|c|c|ll}
\cline{1-3}
\multicolumn{3}{|c|}{2-2-2-1, 2 data}                                                                     &  &  \\ \cline{1-3}
Irreducible components & Loss & \begin{tabular}[c]{@{}c@{}}Is saddle 0 on this \\ component?\end{tabular} &  &  \\ \cline{1-3}
Dim 6, deg 1, 1 component      & 13   & Yes                                                                       &  &  \\ \cline{1-3}
Dim 6, deg 3, 1 component      & 13   & Yes                                                                       &  &  \\ \cline{1-3}
Dim 6, deg 3, 1 component      & 13   & Yes                                                                       &  &  \\ \cline{1-3}
\end{tabular}
\end{center}
\vskip -0.1in
\end{minipage}\hfill
\begin{minipage}[t]{\columnwidth}
\caption{The saddle components of 2-2-2 linear net, with 2 training data under condition $\text{rank}(W_2 W_1) < 2$}
\label{table:saddle-2-2-2-critical-2data}
\vskip 0.15in
\begin{center}
\begin{tabular}{|c|c|c|ll}
\cline{1-3}
\multicolumn{3}{|c|}{2-2-2, 2 data}                                                                       &  &  \\ \cline{1-3}
irreducible components & loss & \begin{tabular}[c]{@{}c@{}}is saddle 0 on this \\ component?\end{tabular} &  &  \\ \cline{1-3}

dim3, deg2,1 comp      & 5.3  & No                                                                        &  &  \\ \cline{1-3}
dim3, deg2,1 comp      & 1.69 & No                                                                        &  &  \\ \cline{1-3}
dim0, deg1, 1comp      & 7    & Yes                                                                       &  &  \\ \cline{1-3}
\end{tabular}
\end{center}
\vskip -0.1in
\end{minipage}

\caption{The eigenvalues of a generic point on global minimum  component, one of saddle components and saddle point $\boldsymbol{0}$} 
\label{table:1-1-1-1-1-eigen}
\vskip 0.15in
\begin{center}
\begin{tabular}{|c|c|l|l}
\cline{1-3}
\multicolumn{3}{|c|}{1-1-1-1-1, 2 data}                                    &  \\ \cline{1-3}
Dim 3, deg 4, 1 component                                                                                                                           & Dim 2, deg 1, 6 components                                                                                                                             & Saddle point $\boldsymbol{0}$                                                             &  \\ \cline{1-3}
\begin{tabular}[c]{@{}c@{}}-1.9149e+01 - 4.4317e+01i\\ 6.3600e-07 + 2.4911e-06i\\ -2.7642e-07 + 1.6268e-07i\\ -5.0118e-07 - 2.7757e-08i\end{tabular} & \begin{tabular}[c]{@{}c@{}}-5.0525e+01 - 7.6442e+01i\\ 1.0105e+01 + 1.5288e+01i\\ 3.8473e-33 - 3.2474e-33i\\ -7.6581e-34 - 6.0445e-34i\end{tabular} & \multicolumn{1}{c|}{\begin{tabular}[c]{@{}c@{}}0\\ 0\\ 0\\ 0\end{tabular}} &  \\ \cline{1-3}
\end{tabular}
\end{center}
\vskip -0.1in

  \vspace{1cm} 

\caption{The eigenvalues of a generic point on global minimum component, one of saddle components and saddle point $\boldsymbol{0}$}
\label{table:2-1-2-eigen-2data}
\vskip 0.15in
\begin{center}
\begin{tabular}{|c|c|c|}
\hline
\multicolumn{3}{|c|}{2-1-2, 2data}                                        \\ \hline
Dim 1, deg 1, 2 components                                                                                                                                       & Dim 1, deg 2, 1 component                                                                                                                                        & Saddle point $\boldsymbol{0}$                                                                     \\ \hline
\begin{tabular}[c]{@{}c@{}}9.0735e+00 + 1.6091e-01i\\   -9.0345e+00 + 1.6022e-01i\\   -2.3606e-16 + 4.7020e-17i\\    3.8988e-02 + 3.2113e-01i\end{tabular} & \begin{tabular}[c]{@{}c@{}}-1.8925e+01 + 2.4843e+00i\\   -3.5563e-01 + 1.0326e-01i\\   -1.0784e-06 - 2.4537e-07i\\   -5.6025e+00 - 1.7950e+00i\end{tabular} & \begin{tabular}[c]{@{}c@{}}-9.0554\\   -0.0000\\         0\\    9.0554\end{tabular} \\ \hline
\end{tabular}
\end{center}
\vskip -0.1in

  \vspace{1cm} 
\caption{The eigenvalues of a generic point on global minimum component, one of saddle components and saddle point $\boldsymbol{0}$}
\label{table:2-2-2-eigen-2data}
\vskip 0.15in
\begin{center}
\begin{tabular}{|c|c|c|}
\hline
\multicolumn{3}{|c|}{2-2-2, 2 data}                                       \\ \hline
Dim 4, deg 8, 1 component                                                                                                                                                                                                                                                                 & Dim 3, deg 2, 2 components                                                                                                                                                                                                                                                              & Saddle point $\boldsymbol{0}$                                                                                                                      \\ \hline
\begin{tabular}[c]{@{}c@{}}-2.6274e+00 - 1.6162e+01i\\   -7.5948e+00 - 8.0849e+00i\\    5.1586e+00 - 5.7933e+00i\\    1.9124e-01 + 2.2839e+00i\\    1.6486e-06 + 1.8576e-06i\\   -4.1784e-07 + 6.3060e-07i\\   -2.2429e-07 - 2.2350e-07i\\   -8.0787e-07 - 1.5855e-06i\end{tabular} & \begin{tabular}[c]{@{}c@{}}-9.5240e+00 + 1.4045e+00i\\   -8.8515e+00 - 7.4800e+00i\\   -1.7240e-01 - 8.4029e+00i\\    1.3592e+00 - 9.8640e-01i\\   -5.1422e-01 + 5.0474e-01i\\   -4.7988e-07 + 4.8486e-07i\\    4.2581e-08 - 5.4055e-07i\\    4.7988e-07 - 4.8486e-07i\end{tabular} & \begin{tabular}[c]{@{}c@{}}-4.6056\\   -4.6056\\   -2.6056\\   -2.6056\\    2.6056\\    2.6056\\    4.6056\\    4.6056\end{tabular} \\ \hline
\end{tabular}
\end{center}

\bigskip

\bigskip

\caption{The linear neural variety of 2-2-2 linear residual net, with 2 training data}
\label{table:2-2-2-residue-critical-2data}
\vskip 0.15in
\begin{tabular}{|c|c|c|ll}
\cline{1-3}
\multicolumn{3}{|c|}{2-2-2, 2 data}                                                                       &  &  \\ \cline{1-3}
Irreducible components & Loss & \begin{tabular}[c]{@{}c@{}}Is saddle $\boldsymbol{0}$ on this \\ component?\end{tabular} &  &  \\ \cline{1-3}
Dim 4, deg 8, 1 component      & 0    & No                                                                        &  &  \\ \cline{1-3}
Dim 3, deg 2, 1 component      & 5.3  & No                                                                        &  &  \\ \cline{1-3}
Dim 3, deg 2, 1 component      & 1.69 & No                                                                        &  &  \\ \cline{1-3}
Dim 0, deg 1, 1 component      & 7    & Yes                                                                       &  &  \\ \cline{1-3}
\end{tabular}

\vskip -0.1in

\end{table*}

\subsection*{Critical points arrangement}
When the width of linear network is one and $k=\text{min} \{ d_x, d_y \}$, there are only global minimum component and saddle components passing through point $\boldsymbol{0}$. If we calculate rank zero saddle points, the result is saddle components passing through point $\boldsymbol{0}$. For example, for net with structure "2-2-2-1, 2 training data", after we calculate gradient formula accompany by saddle condition $\text{rank}(W_3 W_2 W_1) < 1$, the results is in Table \ref{table:saddle-2-2-2-1-critical-2data}. Component membership testing is used to determine whether these saddle components equal to saddle components obtained before. 

When the width of linear network is greater than one and $k=\text{min} \{ d_x, d_y \}$, we take net with structure "2-2-2, 2 training data" as example. The saddle points under rank condition $\text{rank}(W_2 W_1) < 2$ are in Table \ref{table:saddle-2-2-2-critical-2data} and the saddle points under rank condition $\text{rank}(W_2 W_1) < 1$ are just components passing through $\boldsymbol{0}$. It can be seen that these results are consistent with Hypothesis \ref{conj:1}.

\subsection*{The number of constant zero eigenvalues on a component}
After we have all the linear neural varieties of these nets, we can calculate the eigenvalues of a sampled generic point. We list three typical results in Table \ref{table:1-1-1-1-1-eigen} \ref{table:2-1-2-eigen-2data} \ref{table:2-2-2-eigen-2data}. In all three cases, the number of constant zero eigenvalues of a generic point is identical with the dimension of component. Table \ref{table:1-1-1-1-1-eigen} shows that saddle point $\boldsymbol{0}$ is totally flat (all eigenvalues are zero). Table \ref{table:2-1-2-eigen-2data} is a example that saddle point $\boldsymbol{0}$ is strict saddle but has more zero eigenvalues than any other point of component and Table \ref{table:2-2-2-eigen-2data} demonstrates that point $\boldsymbol{0}$ is a generic point if it is not a intersection point.

\subsection*{The case of linear residual network}
Due to the landscape translation relation of linear residual network and linear network, the solution of the linear neural variety of linear residual network has components with the same dimensions and degrees as corresponding counterpart. This is shown in Table \ref{table:2-2-2-residue-critical-2data}.

\section{Conclusion and discussion}

In this paper, we use numerical algebraic geometry to get some complex solutions to gradient formulas for multi-layer linear networks although the actual critical points are over the real number field. Three hypotheses are putting forward.
The reasonableness of these three hypotheses is verified by numerical algebraic geometry computation.
Hypothesis \ref{conj:eigenvalue} means that each irreducible saddle component is regarded as a high-dimensional generalization of an isolated saddle point. In addition, every point of a dense subset of each component has at least one strictly negative eigenvalue, that is to say, almost every saddle point is not a problem with SGD by work of \cite{lee2016gradient}. Its Corollary \ref{cor:2} means that the loss surface of most linear networks will be flat near the saddle point 0, and the residual connection will make the path of SGD avoid passing through this area. Hypothesis \ref{conj:constant loss} makes it possible to define the loss of each individual component. Hypothesis \ref{conj:1} implies that the losses of all irreducible components are bounded by the losses of saddle point $0$ and global minimum component. According these three hypotheses we can think of loss surface as a concave surface, and the existence of saddle point $\boldsymbol{0}$ leads to the bottom bulge as Figure \ref{fig:1}. All other components are in the bottom bulge area, which represents that their losses have upper and lower bound.


\nocite{langley00}

\bibliography{icdp2009}

\end{document}